\documentclass[11pt]{amsart}
\usepackage{latexsym,amssymb,amsfonts,amsmath,graphicx,fullpage,url}
\usepackage[vcentermath,enableskew]{youngtab}
\usepackage{color}
\usepackage{cite}
\usepackage[all]{xy}
\usepackage{verbatim}
\newtheorem{theorem}{Theorem}[section]
\newtheorem{proposition}{Proposition}[section]
\newtheorem{lemma}[theorem]{Lemma}
\newtheorem{corollary}[theorem]{Corollary}
\theoremstyle{definition}
\newtheorem{definition}[theorem]{Definition}

\theoremstyle{remark}

\numberwithin{equation}{section}
\def\F{\mathcal{F}}

\def\F{{\mathcal{F}}}
\def\C{{\mathcal{C}}}
\def\S{{\mathcal{S}}}

\def\V{{\mathcal{V}}}

\def\t{{\mathcal{S}(\beta)}}
\def\P{{\mathcal{P}}}
\def\ini{{\rm in }}
\def\fin{{\rm fin}}
\def\R{\mathbb{R}}

\newcommand{\be}{\begin{equation}}
\newcommand{\ee}{\end{equation}}

\newcommand{\bd}{\begin{definition}}
\newcommand{\ed}{\end{definition}}
\newcommand{\bt}{\begin{theorem}}
\newcommand{\et}{\end{theorem}}
\newcommand{\bl}{\begin{lemma}}
\newcommand{\el}{\end{lemma}}
\newcommand{\bp}{\begin{proposition}}
\newcommand{\ep}{\end{proposition}}
\newcommand{\bc}{\begin{corollary}}
\newcommand{\ec}{\end{corollary}}

\newtheorem*{theorem1*}{Theorem \ref{thm:main}}

\newcommand{\old}[1]{}

\author{Karola M\'esz\'aros}
\address{Karola M\'esz\'aros, Department of Mathematics, Cornell University, Ithaca NY 14853  \newline karola@math.cornell.edu
}
\thanks{The author was partially supported by  a National Science
  Foundation Grant  (DMS 1501059).}

\title{Pipe dream complexes and triangulations of root polytopes  belong together}

\begin{document}

\begin{abstract} We show that the pipe dream complex associated to the permutation  $1\text{ } n \text{ }n-1\text{ } \cdots \text{ }2$ can be geometrically realized as a triangulation of the vertex figure of a root polytope. Leading up to this result we show that the Grothendieck polynomial specializes to the $h$-polynomial of the corresponding pipe dream complex, which in certain cases equals the $h$-polynomial of canonical triangulations of root (and flow) polytopes, which in turn equals a specialization of the reduced form of a monomial in the subdivision algebra of root (and flow) polytopes. Thus, we connect Grothendieck polynomials to reduced forms in subdivision algebras and root (and flow) polytopes. We also show that  root polytopes can be seen as projections of  flow polytopes, explaining that these families of polytopes possess the same subdivision algebra.   
\end{abstract}

 \maketitle
 \tableofcontents 

\section{Introduction}
\label{sec:intro}
In this paper we journey from Grothendieck polynomials to  geometric realizations of pipe dream complexes via root polytopes. On this journey we meet reduced forms of monomials in the subdivision algebra of root and flow polytopes, and root and flow polytopes themselves. While the connection between Grothendieck polynomials and pipe dream complexes is a well known one, the other objects in the above list are not universally thought of as tied to 
Grothendieck polynomials and pipe dream complexes. As this work will illustrate, they might indeed belong together.

Grothendieck polynomials  represent K-theory classes on the flag manifold; they generalize Schubert polynomials, which in turn generalize Schur polynomials. We show that Grothedieck polynomials specialize to  $h$-polynomials of pipe dream complexes. Since pipe dream complexes are known to be homeomorphic to balls (except in a trivial case), we get that  their $h$-polynomials, and thus shifted specialized   Grothedieck polynomials  have nonnegative coefficients. Such property was first observed by Kirillov \cite{k2}, who indicated that he had an algebraic proof in mind. 

In \cite{k2} Kirillov also observed that a certain specialization of the shifted Grothendieck polynomial equals a specialization of a particular reduced form in the subdivision algebra of root and flow polytopes. His observation was based on numerical evidence. We explain this equality in terms of the geometry of the underlying pipe dream complex and root (and flow) polytopes. Indeed, we show that the mentioned pipe dream complex can be realized as  the canonical triangulation of the vertex figure of the root polytope. No wonder then the specialized Grothendieck polynomial and reduced form are equal: they are the $h$-polynomial of the pipe dream complex and the $h$-polynomial of the canonical triangulation of the vertex figure of  the root polytope, respectively. That the reduced form can be seen as the $h$-polynomial of the canonical triangulation of the flow polytope, and thus of the canonical triangulation of the vertex figure of the root polytope, was proved in \cite{h-poly1}. The paper \cite{h-poly2} also contains closely related results. 

The outline of the paper is as follows. In Section \ref{sec:groth} we show that  the shifted  $\beta$-Grothendieck polynomial corresponding to the permutation $w$ is the $h$-polynomial of pipe dream complex of $w$ denoted by $PD(w)$. In Section \ref{sec:red} we define the subdivision algebra and reduced forms and recall related results.  We also allude to the connection of Grothendieck polynomials and reduced forms. In Section \ref{sec:r-f} we explain why the subdivision algebras of root and flow polytopes are the same, by showing that  root polytopes are projections of flow polytopes. Finally, in Section \ref{sec:pipe} we tie all the above together, by showing that  the pipe dream complex $PD({1\text{ } n \text{ }n-1\text{ } \cdots \text{ }2})$ can be realized as a canonical triangulation of a vertex figure of a root polytope.  $PD({1\text{ } n \text{ }n-1\text{ } \cdots \text{ }2})$  has been realized previously via the classical associahedron  \cite{assoc, cesar, subwordcluster}.

\section{Grothendieck polynomials}
\label{sec:groth}

In this section we define Grothendieck  polynomials and explain that they  specialize to $h$-polynomials of certain simplicial complexes called pipe dream complexes.   Since the pipe dream complex is homeomorphic to a ball, its $h$-polynomial has nonnegative coefficients. Therefore, we immediately obtain nonnegativity properties of Grothendieck polynomials, which were observed by Kirillov  in \cite{k2}. 

\begin{figure}
\begin{center}
\includegraphics[scale=.65]{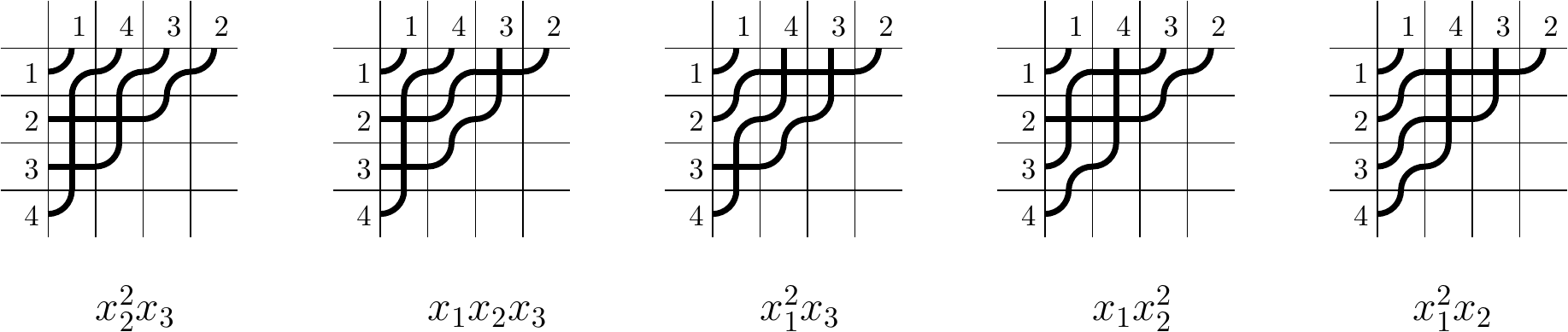}
 \caption{All reduced pipe dreams for $w=1432$ (that is pipe dreams with exactly $3$ crosses). The weights $wt_{x,y}(P)$ when $\bf{y=0}$ are written below the reduced pipe dreams. }
 \label{pipe}
 \end{center}
\end{figure}

There are several ways to express  Grothendieck  polynomials, and we will present the expression in terms of pipe dreams here. 
Given a permutation $w$ in the  symmetric group $S_n$ it can be represented by a triangular table filled with \includegraphics[scale=.5]{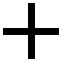}'s (crosses)
 and \includegraphics[scale=.5]{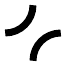}'s (elbows)
 such that (1) the pipes intertwine according to $w$ and (2)  two pipes cross at most once.  Such representations of $w$ are called {\bf reduced pipe dreams}, see Figure \ref{pipe}. Pipe dreams are also known as RC-graphs, and  the reduced pipe dreams of a permutation were shown to be  connected by ladder and chute moves  by Bergeron and Billey in \cite{rc}. To each reduced pipe dream we can associate the {\bf weight} $wt_{x,y}(P):=\prod_{(i,j) \in {\rm cross}(P)} (x_i-y_j)$, with ${\rm cross}(P)$ being the set of positions where $P$ has a cross.  Note that throughout the literature the definition of the weight varies; however, all results can be phrased using any one convention.

\begin{figure}
\begin{center}
\includegraphics[scale=.65]{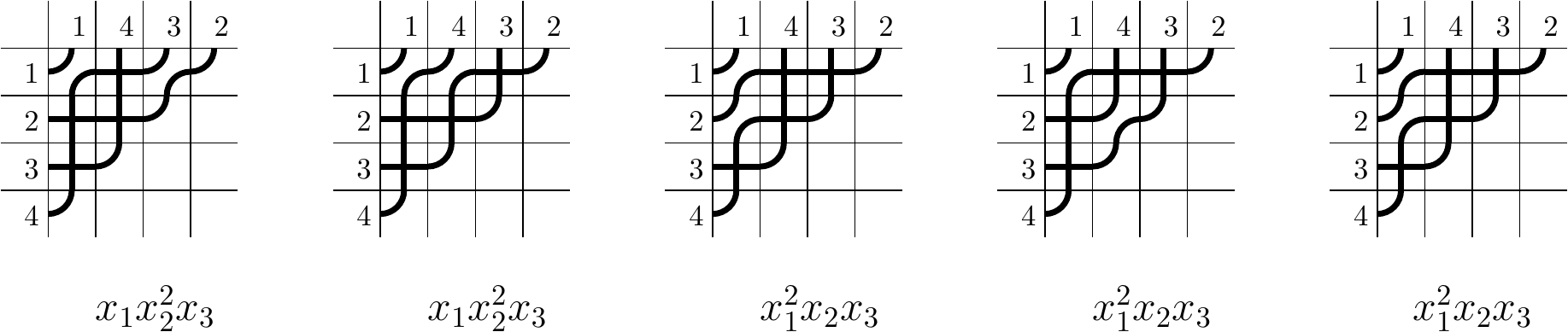}
 \caption{All  pipe dreams with exactly $4$ crosses for $w=1432$. The weights $wt_{x,y}(P)$ when $\bf{y=0}$ are written below the pipe dreams.}
 \label{nonred-pipe}
 \end{center}
\end{figure}

\begin{figure}
\begin{center}
\includegraphics[scale=.65]{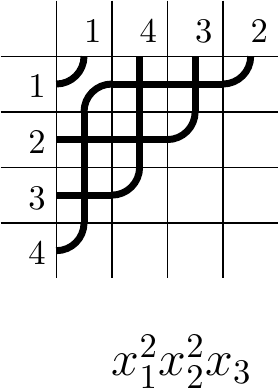}
 \caption{The unique pipe dreams with exactly $5$ crosses for $w=1432$. The weights $wt_{x,y}(P)$ when $\bf{y=0}$ are written below the pipe dream.}
 \label{nonred-one}
 \end{center}
\end{figure}

A  \textbf{nonreduced pipe dream} for $w \in S_n$ is  a triangular table filled with crosses
 and elbows so that (1) the pipes intertwine according to $w$ whereby if two pipes have already crossed previously then we simply ignore the extra crossings and (2) there are two pipes that cross at least twice. All nonreduced pipe dreams for $1432$ with a total of $4$ crosses can be seen  in  Figure \ref{nonred-pipe} and the unique pipe dreams for $1432$ with a total of $5$ crosses can be seen  in  Figure \ref{nonred-one}. We associate a  weight $wt_{x,y}(P)$ to each pipe dream as above.

\begin{figure}
\begin{center}
\includegraphics[scale=.3]{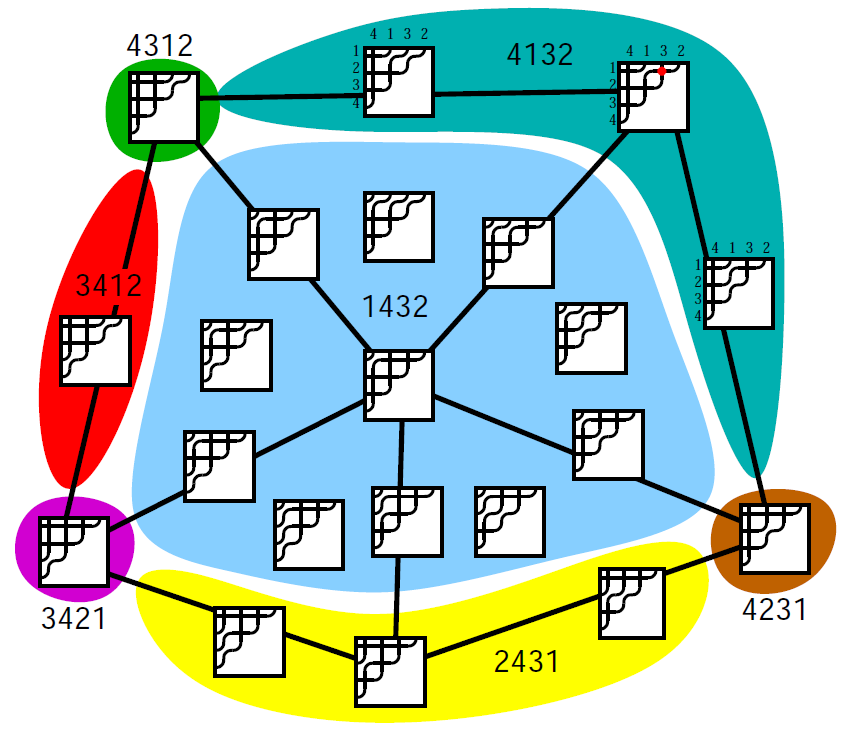}
 \caption{The pipe dream complex $PD({1432})$. Figure used with permission from \cite{allen}.}
 \label{fig:1432}
 \end{center}
\end{figure}

The set ${{\rm Pipes}(w)}$, which  is the set of all pipe dreams of $w$ (both reduced and nonreduced), naturally labels the interior simplices of the {pipe dream complex} $PD(w)$ associated to a permutation $w\in S_n$; see Figure \ref{fig:1432} for $PD({1432})$.  The pipe dream complex $PD(w)$ is a special case of a subword complex and can be defined as follows. A word of size $m$ is an ordered sequence $Q=(\sigma_1, \ldots, \sigma_m)$ of elements from the simple reflections $\{s_1, \ldots, s_{n-1}\}$ in $S_n$. An ordered subsequence $R$ of $Q$ is called a subword of $Q$. A subword $R$ of $Q$ represents $w \in S_n$ if the ordered product of simple reflections in $R$ is a reduced decomposition for $w$. The word $R$ contains $w \in S_n$ if some subsequence of $R$ represents $w$. The  \textbf{pipe dream complex} $PD(w)$ is the set of subwords $(s_{n-1}, s_{n-2}, \ldots, s_1,s_{n-1}, s_{n-2}, \ldots, s_2, s_{n-1}, s_{n-2}, \ldots, s_3, \ldots, s_{n-1}, s_{n-2},  s_{n-1})\backslash R$ whose complements $R$ contain $w$. The word $Q=(s_{n-1}, s_{n-2}, \ldots, s_1,s_{n-1}, s_{n-2}, \ldots, s_2, s_{n-1}, s_{n-2}, \ldots, s_3, \ldots, $ $s_{n-1}, s_{n-2},  s_{n-1})$ is called the triangular word.

The following theorem provides a combinatorial way of thinking about double  Grothendieck polynomials.

\bt \label{allen} \cite{subword, fom-kir} The \textbf{double  Grothendieck polynomial $\mathfrak{G}_w({\bf x, y})$} for $w \in S_n$, where ${\bf x}=(x_1, \ldots, x_{n-1})$ and ${\bf y}=(y_1, \ldots, y_{n-1})$ can be written as

\begin{equation} \mathfrak{G}_w({\bf x, y})=\sum_{P\in {\rm Pipes}(w)}(-1)^{codim_{PD(w)}F(P)} wt_{x,y}(P),
\label{eq:groth}
\end{equation}

\noindent where  ${\rm Pipes}(w)$ is the set of all pipe dreams of $w$ (both reduced and nonreduced), $F(P)$ is the interior face in $PD(w)$ labeled by the pipe dream $P$,  $codim_{PD(w)} F(P)$ denotes the codimension of $F(P)$ in $PD(w)$ and  $wt_{x,y}(P)=\prod_{(i,j) \in {\rm cross}(P)} (x_i-y_j)$, with ${\rm cross}(P)$ being the set of positions where $P$ has a cross. \et

In the spirit of Theorem \ref{allen}, we use the following definition for the {\bf double  $\beta$-Grothendieck polynomial}:

\begin{equation} \mathfrak{G}^{\beta}_w({\bf x, y})=\sum_{P\in {\rm Pipes}(w)}\beta^{codim_{PD(w)}F(P)} wt_{x,y}(P).
\label{eq:bgroth}
\end{equation}

Note that if we assume that $\beta$ has degree $-1$, while all other variables are of degree $1$, then the powers of $\beta$'s simply  make the polynomial $ \mathfrak{G}_w^{\beta}({\bf x, y})$ homogeneous. We chose this definition of $\beta$-Grothendieck polynomials, as it will be the most convenient notationwise for our purposes.

 Next we state a special case of \eqref{eq:bgroth}, since it will play a special role in this section.

\begin{lemma} \label{qt} Denoting  $\mathfrak{G}_w^{\beta}({\bf x, y})$ by $\mathfrak{G}_w^{\beta}({q, t})$  when we set all components of ${\bf x}$ to $q$ and all components of ${\bf y}$ to $t$, we have \begin{equation} \mathfrak{G}_w^{\beta}({q, t})=(q-t)^{l(w)}\sum_{P\in {\rm Pipes}(w)} [\beta(q-t)]^{codim_{PD(w)}F(P)},
\label{eq:g}
\end{equation} 
where $l(w)$ is the length of the permutation $w$, $F(P)$ is the interior face in $PD(w)$ labeled by the pipe dream $P$ and $codim_{PD(w)} F(P)$ denotes the codimension of $F(P)$ in $PD(w)$. 
\end{lemma}

\proof By \eqref{eq:bgroth} we have that

 \begin{equation} \mathfrak{G}_w^{\beta}({q, t})=\sum_{P\in {\rm Pipes}(w)}\beta^{codim_{PD(w)}F(P)} (q-t)^{|{\rm cross}(P)|}.
\end{equation}
Since  the number of crosses in a pipe dream $P$ is  $l(w)+codim_{PD(w)}F(P)$, equation \eqref{eq:g} follows. \qed

The next lemma 
follows from the well-known relation between $f$- and $h$-polynomials. We note that we take $h(\C, x)=\sum_{i=0}^d h_{i}x^{i}$ to be the  $h$-polynomial of a $(d-1)$-dimensional simplicial complex $\C$. For a direct proof of the lemma see \cite{h-poly1}.

\begin{lemma} \cite{Stcom} \label{pure} Let $\C$ be a $(d-1)$-dimensional  simplicial complex homeomorphic to a ball and $f_i^\circ$ be the number of interior faces of $\C$ of dimension $i$. Then 
\begin{equation} \label{h} h(\C, \beta+1)=\sum_{i=0}^{d-1} f_i^\circ \beta^{d-1-i}
\end{equation}
\end{lemma}

Using that the interior simplices of $PD(w)$ are in bijection with pipe dreams of $w$ we obtain the following corollary of Lemma \ref{pure}. 

\begin{corollary} \label{cor}Given $w \in S_n$ we have 

\begin{equation} \label{h1} h(PD(w), \beta+1)=\sum_{P\in {\rm Pipes}(w)}  \beta^{codim_{PD(w)} F(P)},
\end{equation}
where $F(P)$ is the interior face in $PD(w)$ labeled by the pipe dream $P$ and $codim_{PD(w)} F(P)$ denotes the codimension of $F(P)$ in $PD(w)$. 
\end{corollary}

Finally, we obtain the following as a corollary of the above.

\begin{theorem} \label{thm:g-h} We have

\begin{equation} \mathfrak{G}_w^{\beta-1}({q, q-1})= h(PD(w), \beta),
\label{eq:groth-h}
\end{equation}
where $h(PD(w), \beta)$ is the $h$-polynomial of $PD(w)$. In particular we have that $ \mathfrak{G}_w^{\beta-1}({q, q-1})\in \mathbb{Z}_{\geq 0}[\beta]$. 
\end{theorem}

\proof
Corollary \ref{cor} yields $h(PD(w), \beta+1)=\sum_{P\in {\rm Pipes}(w)}  \beta^{codim_{PD(w)} F(P)}$. Together with Lemma \ref{qt} applied when $t=q-1$ we get \eqref{eq:groth-h} in Theorem \ref{thm:g-h}. The nonnegativity of the coefficients of  $ \mathfrak{G}_w^{\beta-1}({q, q-1})$ then follows because of the nonnegativity of the $h$-polynomial of a simplicial complex which is homeomorphic to a ball. Recall that $PD(w)$ is known to be homeomorhpic to a ball, except in the trivial case when it is a $(-1)$-sphere, which case can be checked separately. 
\qed

The nonnegativity of the coefficients of   $\mathfrak{G}_w^{\beta-1}({1, 0})$   was  observed by Kirillov \cite{k2}. Equation \eqref{eq:groth-h} makes clear why this is the case: because it is the $h$-polynomial of a simplicial complex which is homeomorphic to a ball, implying that its coefficients are nonnegative \cite{Stcom}.

\section{Reduced forms in the subdivision algebra}
\label{sec:red}
In this section we point to a connection between reduced forms in the so called subdivision algebra and Grothendieck polynomials.  In  Section \ref{sec:pipe} we provide a  geometric realization of the pipe dream complex $PD({1 n (n-1)\ldots 2})$ via a triangulation of a  root (or flow) polytope, which implies this connection. In Section \ref{sec:r-f} we explain the connection between root and flow polytopes via the subdivision algebra also explaining the algebra's name.

The \textbf{subdivision algebra} $\t$ is a commutative algebra generated by the variables $x_{ij}$, $1\leq i<j\leq n$,  over $\mathbb{Q}[\beta]$, subject to the relations $x_{ij} x_{jk}=x_{ik}(x_{ij}+x_{jk}+\beta)$, for $1\leq i<j<k\leq n$.  This algebra is called the subdivision algebra,  because its relations can be seen geometrically  as subdividing flow  and root polytopes. This is explained in detail in Section \ref{sec:r-f}. The subdivision algebra has been used extensively for subdividing root and flow polytopes  in \cite{prod, mm, h-poly1, h-poly2, root1, root2}.

 A  \textbf{reduced form} of the monomial   in the algebra $\t$  is a polynomial   obtained by successively substituting $x_{ik}(x_{ij}+x_{jk}+\beta)$ in place of an occurrence of $x_{ij} x_{jk}$ for some $i<j<k$ until no further reduction is possible. Note that the reduced forms are not necessarily unique.
  
A possible sequence of reductions in algebra $\t$ yielding a reduced form of $x_{12}x_{23}x_{34}$ is given by

\begin{eqnarray} \label{ex1}
 x_{12} \mbox {\boldmath$  x_{23}x_{34}$} & \rightarrow & \mbox{\boldmath$x_{12}$}x_{24}\mbox{\boldmath$x_{23}$}+\mbox{\boldmath$x_{12}$}x_{34}\mbox {\boldmath$x_{24}$}+\beta \mbox {\boldmath$x_{12}x_{24}$} \nonumber \\
& \rightarrow& \mbox {\boldmath$ x_{24}$} x_{13}\mbox {\boldmath$x_{12}$}+x_{24}x_{23}x_{13}+  \beta x_{24}x_{13}+x_{34}x_{14}x_{12}+x_{34}x_{24}x_{14} \nonumber \\
& &+\beta x_{34}x_{14}+\beta x_{14}x_{12}+\beta x_{24}x_{14}+\beta^2 x_{14} \nonumber \\
& \rightarrow  &x_{13}x_{14}x_{12}+x_{13}x_{24}x_{14}+\beta x_{13}x_{14}+x_{24}x_{23}x_{13}+\beta x_{24}x_{13}\nonumber \\
& & +x_{34}x_{14}x_{12}+x_{34}x_{24}x_{14}+\beta x_{34}x_{14}+\beta x_{14}x_{12}+\beta x_{24}x_{14}\nonumber \\
& &+\beta^2 x_{14}
\end{eqnarray}

\noindent where the pair of variables on which the reductions are performed is in boldface. The reductions are performed on each monomial separately.

 Given a graph $G$, denote by $Q_G(\beta)$ the reduced form of the monomial $\prod_{(i,j) \in E(G)}x_{ij}$ specialized at $x_{ij}=1$ for all  $1\leq i<j\leq n$. The polynomial $Q_G(\beta)$ is unique, though the reduced form with variables $x_{ij}$ is not \cite{root1}.  In recent work \cite{h-poly1} the author  connected $Q_G(\beta)$ to the $h$-polynomials of triangulations of flow polytopes of $\tilde{G}=(V(G)\cup \{s,t\}, E(G)\cup \{(s,i), (i,t) \mid i \in V(G)\})$. Flow polytopes are defined in Section \ref{subsec:fp}; in the next theorem we treat their triangulations, which we denote by $\C$,  as a simplicial complex. Since $\C$ is a simplicial complex homeomorphic to a ball, it follows that  the $h$-polynomial of a triangulation of a polytope  has nonnegative coefficients \cite{Stcom}.

\begin{theorem} \cite{h-poly1} \label{h2} For any graph $G$ we have \be \label{Q} Q_G(\beta)=h(\mathcal{C}, \beta+1),\ee where $\mathcal{C}$ is any unimodular triangulation of the flow polytope $\F_{\tilde{G}}(1,0,\ldots, 0,-1)$ and $h(\mathcal{C}, x)$ is its $h$-polynomial.  In particular, the reduced form $Q_G(\beta-1)$ is a polynomial in $\beta$ with nonnegative coefficients. 
\end{theorem}


For brevity, use the notation $\mathfrak{G}_w(\beta)$  for $\mathfrak{G}_w^{\beta}({1, 0})$, the double $\beta$-Grothendieck polynomial  evaluated when all $x$'s are set to $1$ and $y$'s are set to $0$.  In this notation Theorem \ref{thm:g-h} specialized at $q=1$ states that $\mathfrak{G}_w(\beta)=h(PD(w), \beta+1)$. 

Note the similarity of the statements of Theorems \ref{h2} and \ref{thm:g-h} as a certain polynomial equaling the $h$-polynomial of a simplicial complex. Paired with Kirillov's observation in \cite[Proposition 3.1]{k2} that \be \label{equal} Q_{P_n}(\beta)= \mathfrak{G}_{\pi}(\beta),\ee for 
the permutation $\pi=1\text{ } n \text{ }n-1\text{ } \cdots \text{ }2$ and path graph $P_n=([n], \{(i, i+1)| i \in [n-1]\})$, we obtain that $h(PD(\pi), \beta)=h(\C, \beta)$, where $\mathcal{C}$ is any unimodular  triangulation  of the flow polytope $\F_{\tilde{P_n}}$. The previous raises the natural question: can $PD({\pi})$ be realized geometrically as  a  triangulation  $\mathcal{C}$ of the flow polytope $\F_{\tilde{P_n}}$?   The answer is almost yes as we explain in the next sections.

\section{On the relation of root and flow polytopes}
\label{sec:r-f}
This section explains the geometric reasons for root and flow polytopes to have the same subdivision algebras and in turn to possess  dissections with identical descriptions via  reduced forms  \cite{prod, mm, h-poly1, h-poly2, root1, root2}.  The simplest reason for the above would be if root and flow polytopes were equivalent. While this is not  the case, the truth does not lie far from it, as we will see.

\subsection{Root polytopes.} In the terminology of \cite{p1}, a root polytope of   type $A_{n}$ is  the convex hull of the origin and some of the points $e_{ij}^-:=e_i-e_j$ for $1\leq i<j \leq n+1$, where $e_i$ denotes the $i^{th}$ coordinate vector in $\mathbb{R}^{n+1}$. A very special root polytope is the full root polytope  $$\mathcal{P}(A_{n}^+)=\textrm{ConvHull}(0,  e_{ij}^- \mid  1\leq i<j \leq n+1),$$  where  $e_{ij}^-=e_i-e_j$. In this paper we restrict ourself  to a  class of root polytopes including $\mathcal{P}(A_{n}^+)$, which have subdivision algebras \cite{root1}. 

Let $G$ be an acyclic   graph on the vertex set $[n+1]$.  Define $$\mathcal{V}_G=\{e_{ij}^- \mid  (i, j) \in E(G), i<j\}, \mbox{ a set of vectors associated to $G$;}$$

 $$\mathcal{C}(G)=\langle \mathcal{V}_G \rangle :=\{\sum_{ e_{ij}^- \in \mathcal{V}_G}c_{ij} e_{ij}^- \mid  c_{ij}\geq 0\}, \mbox{ the cone associated to $G$; and } $$  
  $$\overline{\mathcal{V}}_G=\Phi^+ \cap \mathcal{C}(G), \mbox{ all the positive roots of type $A_n$ contained in $\mathcal{C}(G)$}, $$
   where $\Phi^+=\{e_{ij}^-  \mid1\leq i<j \leq n+1\}$ is the set of         positive roots of type $A_n$.  
    
The root polytope $\mathcal{P}(G)$ associated to the acyclic graph $G$ is  
 
 \begin{equation} \label{eq11} \mathcal{P}(G)=\textrm{ConvHull}(0, e_{ij}^- \mid e_{ij}^- \in \overline{\mathcal{V}}_G)\end{equation} The root polytope $\mathcal{P}(G)$ associated to graph $G$ can also be defined as \begin{equation} \label{eq21} \mathcal{P}(G)=\mathcal{P}(A_n^+) \cap \mathcal{C}(G).\end{equation}  
   
 Note that $\mathcal{P}(A_{n}^+)=\mathcal{P}(P_{n+1})$ for the  path graph  $P_{n+1}$ on the vertex set $[n+1]$.  
 
We can view reduced forms  in the subdivision algebra in terms of graphs, as hinted at in the previous section. 

    The {\bf reduction rule for graphs:} Given   a graph $G_0$ on the vertex set $[n+1]$ and   $(i, j), (j, k) \in E(G_0)$ for some $i<j<k$, let   $G_1, G_2, G_3$ be graphs on the vertex set $[n+1]$ with edge sets
  \begin{eqnarray} \label{graphs}
E(G_1)&=&E(G_0)\backslash \{(j, k)\} \cup \{(i, k)\}, \nonumber \\
E(G_2)&=&E(G_0)\backslash \{(i, j)\} \cup \{(i, k)\},\nonumber \\ 
E(G_3)&=&E(G_0)\backslash \{(i, j), (j, k)\} \cup \{(i, k)\}. 
\end{eqnarray}

    We say that $G_0$ \textbf{reduces} to $G_1, G_2, G_3$ under the reduction rules defined by equations (\ref{graphs}).
 
The reason for the name subdivision algebra is the following key lemma appearing in \cite{root1}:

    \begin{lemma} \cite{root1} \label{reduction_lemma} \textbf{(Reduction Lemma for Root Polytopes)} 
Given   an acyclic  graph $G_0$ with $d$ edges, let  $(i, j), (j, k) \in E(G_0)$ for some $i<j<k$ and $G_1, G_2, G_3$ as described by equations (\ref{graphs}).   Then  
$$\mathcal{P}(G_0)=\mathcal{P}(G_1) \cup \mathcal{P}(G_2)$$   where all polytopes  $\mathcal{P}(G_0), \mathcal{P}(G_1), \mathcal{P}(G_2)$ are   $d$-dimensional and    
$$\mathcal{P}(G_3)=\mathcal{P}(G_1) \cap \mathcal{P}(G_2)  \mbox{   is $(d-1)$-dimensional. } $$
\end{lemma}
      \medskip
      
       What the Reduction Lemma really says is that performing a reduction on an acyclic graph $G_0$  is the same as dissecting the $d$-dimensional  polytope  $\mathcal{P}(G_0)$ into two $d$-dimensional polytopes $\mathcal{P}(G_1)$ and $ \mathcal{P}(G_2)$, whose vertex  sets are  subsets of the vertex set of   $\mathcal{P}(G_0)$, whose interiors are disjoint, whose union is $\mathcal{P}(G_0)$, and whose intersection is a facet of both.  It is clear then that the reduced form can be seen as a dissection of the root polytope into simplices.

\subsection{Flow polytopes.} \label{subsec:fp} Now we define flow polytopes and explain  the analogue of the Reduction Lemma for them. 
 Let $G$ be a loopless graph on the vertex set $[n+1]$, and let $\ini(e)$ denote the smallest (initial) vertex of edge $e$ and $\fin(e)$ the biggest (final) vertex of edge $e$.  Think of fluid flowing on the edges of $G$ from the smaller to the bigger vertices, so that the total fluid volume entering vertex $1$ is one and leaving vertex $n+1$ is one, and there is conservation of fluid at the intermediate vertices. Formally, a \textbf{flow} $f$ of size one on $G$ is a function $f: E \rightarrow \R_{\geq 0}$ from the edge set $E$ of $G$ to the set of nonnegative real numbers such that 
  
  $$1=\sum_{e \in E, \ini(e)=1}f(e)= \sum_{e \in E,  \fin(e)=n+1}f(e),$$
  
  and for $2\leq i\leq n$
  
  $$\sum_{e \in E, \fin(e)=i}f(e)= \sum_{e \in E, \ini(e)=i}f(e).$$
  
  \medskip

  The \textbf{flow polytope} $\F_G$ associated to the graph $G$ is the set of all flows $f: E \rightarrow \R_{\geq 0}$ of size one.   
  
  In this paper we restrict our attention to flow polytopes of certain augmented graphs $\tilde{G}=(V(G)\cup \{s,t\}, E(G)\cup \{(s,i), (i,t) \mid i \in V(G)\})$:

 \begin{lemma} \cite{prod, mm} \label{red} \textbf{(Reduction Lemma for Flow Polytopes)}  Given a graph $G_0$ on the vertex set $[n+1]$ and   $(i, j), (j, k) \in E(G_0)$,   for some $i<j<k$, let $G_1, G_2, G_3$ be as in equations (\ref{graphs}). Then 
 $$\F_{{\tilde{G}_0}}=\F_{{\tilde{G}_1}} \bigcup \F_{{\tilde{G}_2}},$$ where all polytopes $\F_{{\tilde{G}_0}},\F_{{\tilde{G}_1}}, \F_{{\tilde{G}_2}},$ are of the same dimension and $$\F_{{\tilde{G}_3}}=\F_{{\tilde{G}_1}} \cap\F_{{\tilde{G}_2}}  \mbox{   is  one dimension less. } $$
\end{lemma}

\subsection{Are root polytopes and flow polytopes the same?} Given an acyclic graph $G$ Lemmas \ref{reduction_lemma} and \ref{red} imply that we can dissect $\P(G)$ and $\F_{{\tilde{G}}}$ with identical procedures. Are then $\P(G)$ and $\F_{{\tilde{G}}}$  equivalent for  acyclic graphs $G$?

Note that the dimension of $\P(G)$ is $|E(G)|$, while the dimension of $\F_{{\tilde{G}}}$ is $|E(G)|+|V(G)|-1$, so the polytopes cannot be identical. However, we show that  $\F_{{\tilde{G}}}$ can be projected onto an $|E(G)|$-dimensional polytope $\S(G)$ that is equivalent to 
$\P(G)$.  When with the subdivision algebra we are dissecting  $\P(G)$ and $\F_{{\tilde{G}}}$ in identical ways, we get the corresponding (identifiable) induced dissections on $\S(G)$ and $\P(G)$.

Recall the well-known charaterization of the vertices of flow polytopes. 
 
 \bl \label{vertices}\cite[Section 13.1a]{sch} The vertex set of $\F_G$ are the unit flows on increasing paths going from the smallest to the largest vertex of $G$. 
 \el
 
 The polytope $\F_{{\tilde{G}}}$ naturally lives in the space $\mathbb{R}^{|E(\tilde{G})|}$, with the coordinates corresponding to the edges of $\tilde{G}$.  Denoting by $e_{(i,j)}$ the unit coordinate corresponding to the edge $(i,j) \in E(\tilde{G})$, we see that the vectors $e_{(i,j)}$, $(i,j) \in E(G)$, $e_{(s,i)},e_{(i,t)}$, for $i \in [n]$, are an orthonormal basis of $\mathbb{R}^{|E(\tilde{G})|}$. Projecting onto the subspace $W$  of $\mathbb{R}^{|E(\tilde{G})|}$ spanned by $e_{(i,j)}$, $(i,j) \in E(G)$, let  the polytope $\S(G)$ be the image of  $\F_{{\tilde{G}}}$. Denote the mentioned projection by $p$. The vertices of  $\S(G)$ are $0$ and vertices of the form $e_{(i_1, i_2)}+e_{(i_2, i_3)}+\cdots+e_{(i_k, i_{k+1})}$, where $i_1<\cdots<i_{k+1}$, ${(i_1, i_2)}, {(i_2, i_3)},\ldots, {(i_k, i_{k+1})} \in E(G)$. 
 
 Define the map $f:W\rightarrow \mathbb{R}^n$ as follows: $f(e_{(i,j)})=e_i-e_j$, for $(i,j) \in E(G)$, and extend linearly. It follows by definition that the image of $\S(G)$ under $f$ is $\P(G)$. Since for an acyclic graph $G$ the vectors $e_i-e_j$, $(i,j) \in E(G)$, are linearly independent, we get that $f$ is an affine map which is a bijection onto $\P(G)$ when restricted to $\S(G)$. Thus, $\S(G)$ and $\P(G)$ are affinely (and thus combinatorially) equivalent polytopes. 
 
 Let $G_0$ be an acyclic graph, and let $G_1, G_2, G_3$ be as specified \eqref{graphs} .  Then a check shows that $f(p(\F_{{\widetilde{G_i}}}))=\P(G_i)$, for $i \in [3]$ and thus  any dissection of  $\F_{{\widetilde{G_0}}}$ that we obtain by repeated reductions as in  Lemma \ref{red} under the map $f \circ p$ yields a dissection of $\P(G_0)$ obtained by the same sequence of reductions as interpreted in Lemma  \ref{reduction_lemma}. 
 \medskip

  The above considerations prove the  following theorem, which   relates  root and flow polytopes.  The maps $p$ and $f$ are as defined above. 
 
 \bt  \label{root-flow} The root polytope $\P(G)$ is equivalent to $\S(G)$, which is a projection of  $\F_{{\tilde{G}}}$. Indeed, $\P(G)=f(p(\F_{{\tilde{G}}}))$. Moreover, when the reductions \eqref{graphs} are performed on $G$ yielding  dissections $\mathcal{D}_1$ and $\mathcal{D}_2$ of  $\P(G)$  and  $\F_{{\tilde{G}}}$, respectively, then  $\mathcal{D}_1$  is the image of  $\mathcal{D}_2$ under $f\circ p$. \et

  \medskip
 It is in the sense of Theorem \ref{root-flow} that root polytopes and flow polytopes  of acyclic graphs are the same. Since root polytopes are lower dimensional by definition and in this paper we are only concerned with acyclic graphs, namely, the path graph, we will use root polytopes in the rest of the paper.

 \section{Geometric realization of pipe dream complexes via root polytopes}
\label{sec:pipe}
The main theorem of this section is that the canonical triangulation of the vertex figure $\V(P_n)$ of  $\P(P_n)$ at $0$ is a geometric realization of the pipe dream complex $PD({1\text{ } n \text{ }n-1\text{ } \cdots \text{ }2})$.  The vertex figure of a polytope $P$ at vertex $v$ is the intersection of a hyperplane $\mathcal{H}$ with $P$, such that  vertex $v$ is on one side of $\mathcal{H}$  and all the other vertices of $P$ are on the other side of $\mathcal{H}$. See \cite[p.54]{ziegler} for further details. We now explain the canonical triangulation of $\P(P_n)$; there is an analogous triangulation for all root (and flow) polytopes \cite{root1, h-poly2}, but since we are only concerned with $\P(P_n)$ in this section, we restrict our attention to this case. $PD({1\text{ } n \text{ }n-1\text{ } \cdots \text{ }2})$ has previously been been realized via the classical associahedron   \cite{assoc, cesar, subwordcluster}.

Recall that a graph $G$ on the vertex set $[n]$ is said to be {\bf noncrossing} if there are no vertices $i <j<k<l$ such that $(i, k)$ and $(j, l)$ are edges in $G$. 
A graph $G$ on the vertex set $[n]$ is said to be {\bf alternating} if there are no vertices $ i <j<k $ such that $(i, j)$ and $(j, k)$ are edges in $G$.

\begin{theorem}\cite{GGP, root1} \label{path} Let $T_1, \ldots, T_k$ be all the noncrossing alternating spanning trees of $K_n$. Then $\P(T_1), \ldots, \P(T_k)$ are top dimensional simplices in a triangulation of $\P(P_n)$. Moreover, $$\P(T_{i_1})\cap \cdots \cap \P(T_{i_l})=\P(T_{i_1}\cap \cdots \cap T_{i_l}),$$ where $i_1,\ldots,i_l \in [k]$,  and   $T_{i_1}\cap \cdots \cap T_{i_l}=([n], \{(i,j)| (i,j) \in E(T_{i_1})\cap \cdots \cap E(T_{i_l})$. 
 \end{theorem}
 
 The triangulation described in Theorem \ref{path} is called the \textbf{canonical triangulation} of $\P(P_n)$. Since all top dimensional simplices in it contain $0$, we see that $\V(P_n)$ has a triangulation indexed by the same noncrossing alternating spanning trees:

\begin{theorem} \label{v-path} Let $T_1, \ldots, T_k$ be all the noncrossing alternating spanning trees of $K_n$. Then $\P(T_1)\cap \V(P_n), \ldots, \P(T_k)\cap \V(P_n)$ are top dimensional simplices in a triangulation of $\V(P_n)$. Moreover, $$\P(T_{i_1})\cap \cdots \cap \P(T_{i_l})\cap \V(P_n)=\P(T_{i_1}\cap \cdots \cap T_{i_l})\cap \V(P_n),$$ where $i_1,\ldots,i_l \in [k]$,  and   $T_{i_1}\cap \cdots \cap T_{i_l}=([n], \{(i,j)| (i,j) \in E(T_{i_1})\cap \cdots \cap E(T_{i_l})$. 
 \end{theorem}
 
 We call the triangulation described in Theorem \ref{v-path} the \textbf{canonical triangulation} of $\V(P_n)$. The following is the main theorem of this section. 
 
 \bt \label{thm:gr} The canonical triangulation of $\V(P_n)$ is a geometric realization of the pipe dream complex $PD({1\text{ } n \text{ }n-1\text{ } \cdots \text{ }2})$. 
 \et

\begin{figure}
\begin{center}
\includegraphics[scale=.65]{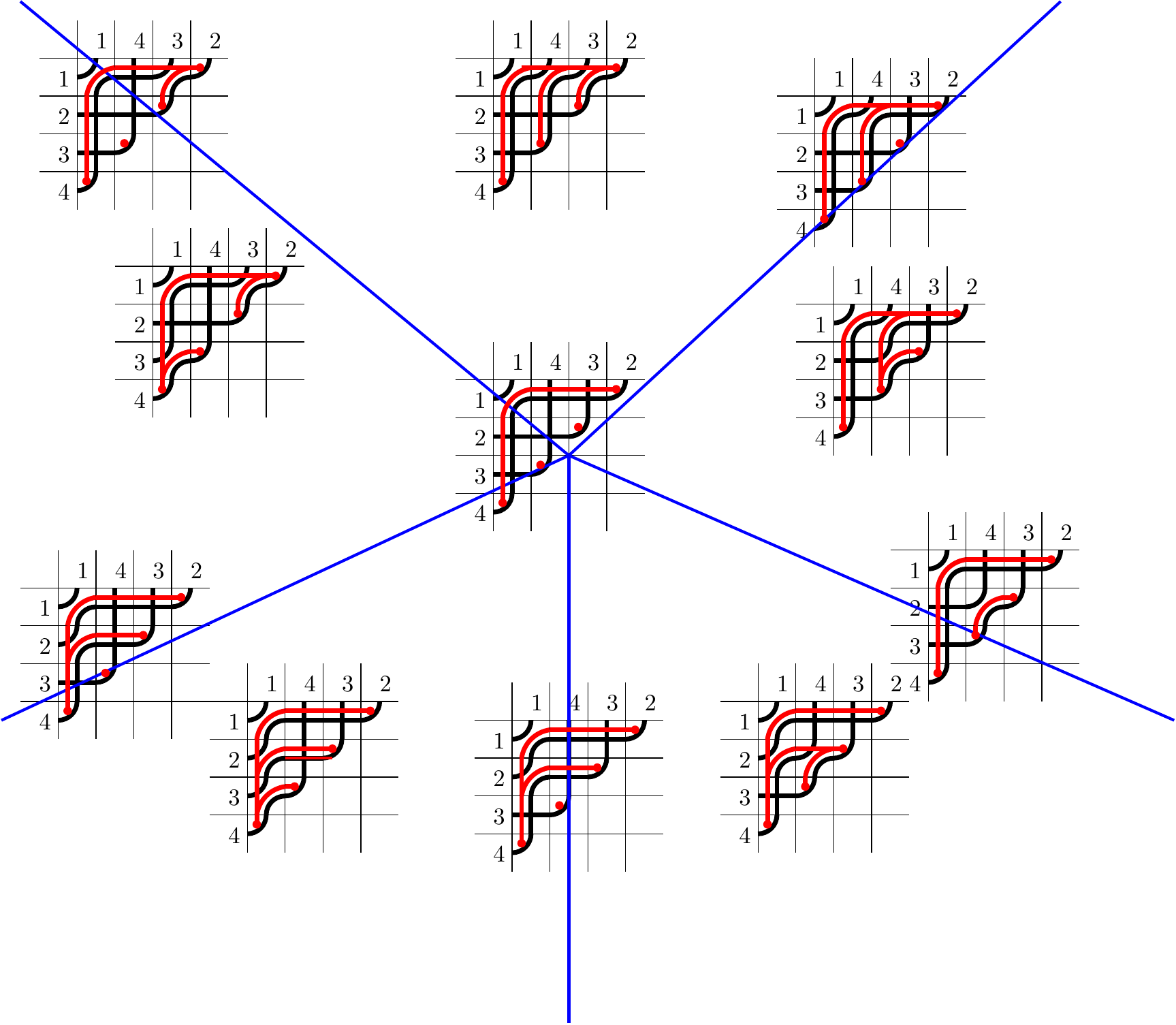}
 \caption{The interior simplices of $PD({1432})$ with the pipe dreams that label them. The graphs obtained via the bijection $G$ can be seen in red  juxtaposed on top of the pipe dreams (rotated by  $45^\circ$). Note that the top dimensional simplices are indeed labeled by the noncrossing alternating spanning trees of $K_n$.}
 \label{bij}
 \end{center}
\end{figure}

Before proceeding to prove Theorem \ref{thm:gr} we note that a proof of it could be obtained using the previous realization of $PD({1\text{ } n \text{ }n-1\text{ } \cdots \text{ }2})$ via the associahedron. However,  instead  we will give a proof using root polytopes. 

\medskip

\noindent \textit{Proof of Theorem \ref{thm:gr}.} First note that the dimensions of both $\V(P_n)$ and $PD({1\text{ } n \text{ }n-1\text{ } \cdots \text{ }2})$ are $n-2$.
Recall that 
 the top dimensional simplices of $PD({1\text{ } n \text{ }n-1\text{ } \cdots \text{ }2})$ are indexed by reduced pipe dreams of ${1\text{ } n \text{ }n-1\text{ } \cdots \text{ }2}$, whereas their intersections by nonreduced pipe dreams in the following way. If we identify a pipe dream $P$ with its set of crosses, then the simplex at the intersection of the simplices labeled by pipe dreams $P_{i_1}, \ldots, P_{i_l}$ is labeled by a pipe dream  $P_{i_1}\cup \cdots \cup P_{i_l}$, where in the latter we simply let the crosses be all the crosses in  $P_{i_1}, \ldots, P_{i_l}$.  See Figure \ref{fig:1432}.
 
 We first show that the top dimensional simplices in   $PD({1\text{ } n \text{ }n-1\text{ } \cdots \text{ }2})$ are in bijection with the top dimensional simplices in the canonical triangulation of $\V(P_n)$. Such a bijective map  $G$ is easy to define. Given a reduced pipe dream $P$, let $$G(P)=([n], \{(i,j) | \mbox{there is an elbow in box } (n-j+1, i) \mbox{ in P}\}).$$ See Figure \ref{bij} for an example. 
 
 The map $G$ is clearly one-to-one. On the other hand we know that both reduced pipe dreams of ${1\text{ } n \text{ }n-1\text{ } \cdots \text{ }2}$ and noncrossing alternating spanning trees of $K_n$ are counted by the Catalan numbers \cite{woo, GGP}, thereby immediately yielding that  $G$ is a bijection between the two sets. 
 
 Given that the simplex at the intersection of the simplices labeled by pipe dreams $P_{i_1}, \ldots, P_{i_l}$ is labeled by a pipe dream  $P_{i_1}\cup \cdots \cup P_{i_l}$ (as explained above),  and  $$\P(T_{i_1})\cap \cdots \cap \P(T_{i_l})\cap \V(P_n)=\P(T_{i_1}\cap \cdots \cap T_{i_l})\cap \V(P_n),$$ where $i_1,\ldots,i_l \in [k]$,  and   $T_{i_1}\cap \cdots \cap T_{i_l}=([n], \{(i,j)| (i,j) \in E(T_{i_1})\cap \cdots \cap E(T_{i_l})$ (as in Theorem \ref{v-path}) we have that the bijection $G$  extends to the lower dimensional interior simplices of $PD({1\text{ } n \text{ }n-1\text{ } \cdots \text{ }2})$  and 
 $\V(P_n)$. See Figure \ref{bij}.  Moreover, the same map also extends to the boundary simplices in the canonical triangulation of $\V(P_n)$ and those in $PD({1\text{ } n \text{ }n-1\text{ } \cdots \text{ }2})$. Therefore, we can conclude that the canonical triangulation of $\V(P_n)$ is a  geometric realization of  $PD({1\text{ } n \text{ }n-1\text{ } \cdots \text{ }2})$. 
  \qed
 
 \medskip 
 
We remark that the $h$-vector of the canonical triangulation of $\P(P_n)$, and so also of $PD({1\text{ } n \text{ }n-1\text{ } \cdots \text{ }2})$ consists of  Narayana numbers; see  \cite[Exercise 6.31b]{ec2}.

\section*{Acknowledgements} I am grateful to Allen Knutson for the many helpful discussions and references about topics related to this research.  Special thanks to Sergey Fomin for an extensive   conversation   about double Grothendieck polynomials. I also thank Lou Billera and Ed Swartz  for several valuable discussions and the anonymous referees for their comments.

\bibliography{biblio-kir}
\bibliographystyle{alpha}

\end{document}